\documentclass[doublespace]{amsart}
\usepackage{graphicx}
\usepackage{amssymb,amscd,amsmath,amsthm}

\newtheorem{theorem}{Theorem}
\newtheorem{corollary}[theorem]{Corollary}
\newtheorem{proposition}[theorem]{Proposition}

\theoremstyle{definition}
\newtheorem{definition}[theorem]{Definition}

\theoremstyle{remark}

\newcommand{\A}{{\Bbb A}}
\newcommand{\C}{{\Bbb C}}

\newcommand{\R}{{\Bbb R}}

\renewcommand{\P}{{\Bbb P}}
\newcommand{\Z}{{\Bbb Z}}

\newcommand{\Res}{\operatorname{Res}}
\newcommand{\Tr}{\operatorname{Tr}}
\newcommand{\rank}{\operatorname{rank}}
\newcommand{\Sol}{\operatorname{Sol}}
\newcommand{\Spec}{\operatorname{Spec}}
\newcommand{\GL}{\operatorname{GL}}

\newcommand{\sC}{\mathcal{C}}

\newcommand{\sO}{\mathcal{O}}
\newcommand{\sV}{\mathcal{V}}

\title[Period determinant of a connection]{Period determinant of an
irregular 
connection over an elliptic curve} 
\author{Byungheup Jun}
\date{12. Feb. 2003} 
\address{School of Mathematics, KIAS, Seoul, Korea}
\email{bhjun@kias.re.kr}

\begin{document}

\begin{abstract}
In this article, we calculate the period of an irregular singular connection 
$\nabla=d+dy$ 
on the legendre curve  $U: y^2 = x(x-1)(x-\lambda)$.
We calculate its de Rham cohomology and the cycles in the homology 
of the dual connection and describe the period matrix.
Terasoma's work is introduced to approximate
the direct image connection $\pi_*(\nabla)$ by a sequence of regular 
connections as an intermediate step where
$\pi: U\to \A^1=\Spec k[y], (x,y)\mapsto y$.  
Finally, we will make the comparison of the period obtained 
by this approximation and that $\nabla$ over $U$.
\end{abstract}
\maketitle
\tableofcontents

\section{Introduction}

In this article, we will calculate the period determinant of the
connection $\nabla=d + dy$ over
the Legendre elliptic curve: $y^2 = x(x-1)(x-\lambda)$. $\nabla$
has an irregular singularity at
$\infty$. The period for irregular connection is defined as a
pairing of de Rham cohomology and
a newly defined homology in \cite{be}. It is shown to be perfect.

For a regular singular connection, the period is defined
as a pairing of de Rham cohomology and the homology
relative to the singularities valued in its dual local system, via a
singular integration.
From the view  of the theory of characteristic classes, it was discussed
profoundly by T.~Saito and T.~Terasoma in \cite{st}.
Especially, on $\P^1$, it is described  explicitly 
as a  product
formula written in terms of the Gamma factors and the tame symbols 
by Terasoma
in \cite{te1}.

As for the irregular connections,  very few cases are known.
On $\A^1- \{ \lambda_1,\ldots,\lambda_n \}$, when the connection
is of the form $d+\sum \frac{s_i}{x-\lambda_i} +dF$ for a polynomial 
$f$ and positive real numbers $s_i$ , the period determinant
is calculated via approximation by series of rank $1$ regular 
singular period also by Terasoma(cf. \cite{te2}).

In our case, the period determinant involves some elliptic integrals,
which are not easily calculated.
The main idea is to calculate the period of its direct image 
connection under the second projection 
$\pi:U\to \A^1,  (x,y)\mapsto y$ in order
to make it a higher rank connection on $\P^1$.
Then at least those integrations
are on $\A^1$, so that
it can be approximated by a
sequence of periods of properly chosen regular
singular connections on $\P^1$ converging to the irregular connection.
Those periods of regular connections
still do not  satisfy some additional conditions for applying 
the product formula: specially
those on the eigenvalues of the residues at the auxiliary 
regular singular points obtained while projecting down. 
Let $D\in \P^1$ be those ramification points and $d_D$ 
be the standard exterior differentiation of $\sO(-D)$ in $\sO$.
Tensored with $d_D$,  $\pi_*(\nabla)$ is approximated by  connections
whose periods are computable using the product formula.

Returning to the original connection, we take another connection
$\nabla_\Sigma= \nabla \otimes d_\Sigma$, where $\Sigma := \pi^{-1} D$.
This new connection on $U$ produces same periods 
as $d_D \otimes \pi_*(\nabla)$ since
they share an isomorphic relative homology group as well as
an isomorphic de Rham cohomology.
We will compare the period of $\nabla$ with 
that of $\nabla_\Sigma$
through the long exact sequences of de Rham cohomology and of the 
homology to produce exact value of the  period determinant of 
$\nabla$ on $U$.
It will be treated in the generic case:$\lambda^2 - \lambda+1 \ne 0$.

At the end,  we  calculate separately the period determinant at the special 
value of $\lambda$.

{\it  Acknowledgement.}
This article is prepared as a part of the Ph.D thesis of the author. 
He thanks to
the thesis advisor H\'el\`ene Esnault for sharing her idea, 
time and  for encouragement.
He also thanks to Marco Hien  and  Spencer
Bloch for showing their interests and useful discussions on this work.

\section[De Rham cohomology and homology]{De Rham cohomology and homology for irregular singular connections}

In this section, we calculate and describe the de Rham cohomology and the homology  
for $\nabla = d + dy$ over an affine legendre elliptic curve $U: y^2 = x(x-1)(x-\lambda)$.
Let $X$ be the completion of $U$ in $\P^2$(ie. $U= X - \{ \infty \}$).
Generally, for an integrable connection, the de Rham cohomology is defined as the 
hypercohomology of the de Rham complex.

Since $\nabla$ admits no solution, $H^0_{dR}(U,\nabla) = 0$ and
\begin{proposition}
$$
H^1_{dR}(U,\nabla) = 
\frac{\Gamma(X,\Omega_X(5\cdot \infty))}{\nabla \Gamma(X,\sO_X(\infty))}
$$
Therefore, $H^1_{dR}(U,\nabla)$ is generated by $\frac{dx}y, x\frac{dx}y, dx, xdx$.
\end{proposition}
\begin{proof}
Since we take the cohomology on an affine variety, we have the identification 
$H^1_{dR}(U, \nabla) = \Gamma(X, \Omega_X(*\infty))/
\nabla\Gamma(X, \sO_X(*\infty))$. It follows that the map
$\Gamma(X,\Omega((n+4)\cdot\infty)) \stackrel{\alpha}{\to} H^1_{dR}(U,\nabla)$ has 
the kernel $\nabla(X,\Omega(n\cdot\infty))$ and is surjective for $n\gg 0$.
As $\dim \alpha \Gamma(X, \Omega((n+4)\cdot\infty)) = 
h^0(X, \Omega((n+4)\cdot\infty)) - h^0(X, \sO(n\cdot\infty)) = 4$ for $n \ge 0$.
Therefore, we have the identification.  
Since $\frac{dx}y, x\frac{dx}y, dx, xdx, dy$ generates $\Gamma(X,\Omega(4\cdot\infty))$ and $dy= \nabla(1)$,
$\frac{dx}y, x\frac{dx}y, dx, xdx$ is a basis for $H^1_{dR}(U,\nabla)$.  
\end{proof}

We denote by $H^{irreg}_*$ the homology for irregular connections defined by Bloch
and Esnault in \cite{be}. This homology is defined and studied over curves. It is known 
to make a perfect pair with the de Rham cohomology of the dual connection. The homology group has cycles that decays rapidly, as well as topological cycles.

From the duality, we see $\dim H^{irreg}_0(U,\nabla^*) = 0$ and 
$\dim H^{irreg}_1(U,\nabla^*) = 4$, where $\nabla^*$ is the dual connection of 
$\nabla=d+dy$ over $U$.

Let $\Delta$ be a small disk around $\infty$. Then the irregular homology  groups
are plugged into the following exact sequence.
\begin{equation}
\begin{split}
0 \to H_1(U,\C) \stackrel{1\otimes \exp y}{\to} H^{irreg}_1(U,\nabla^*) \to
H^{irreg}_1(\Delta^*,\partial \Delta, \nabla^*) \\
\stackrel{\delta}{\to} H_0 (U,\C) \to H^{irreg}_0(U,\nabla^*) \to 0
\end{split}
\end{equation}
In the above, $H_1(U,\C) = H_1(X,\C)$ generated by $\gamma_1, \gamma_2$.
The other cycles in $H^{irreg}_1(U,\C)$ are in the kernel of $\delta$. 
With the identification $H_0(U,\C) = \C <\exp y >$, fixing a local parameter 
$t =-(1/y)^{1/3}$ at
$\infty$, we have the three rapid sectors where the solution decays rapidly:
$$ 
\frac{2(i-1)\pi}3 - \frac\pi6 < \arg t < \frac{2(i-1)\pi}3 + \frac\pi6, \quad\text{for~~~}
i = 1,2,3
$$
We call it the \textit{$i$-th rapid decay sector}. Let $\eta_i$ be a chain from
a fixed point $p$ at $\partial \Delta$ to $\infty$ along the $i$-th rapid decay sector.
$\eta_i \otimes \exp y$ for $i=1,2,3$ generate $H^{irreg}_1(\Delta^*,\partial\Delta,\nabla^*)$.  $\delta$ is the augmentation 
map $\sum_i \alpha_i \eta_i \otimes \exp y \mapsto \sum_i a_i <\exp y>$.
Therefore the kernel is generated by $(\eta_1 -\eta_2)\otimes\exp y$ and
$(\eta_2 - \eta_3)\otimes \exp y$.

Therefore we obtain the following:
\begin{proposition}
Let $\gamma_1, \gamma_2$ be the two generators of $H_1(X,\Z)$ and $\gamma_3$
(resp. $\gamma_4$) be the chain $\eta_2 - \eta_1$(resp. $\eta_3 - \eta_2$). 
$H^{irreg}_1(U,\nabla^*)$ has a basis
$ \gamma_i \otimes \exp y$ for $i=1,2,3,4$.
\end{proposition}

Let $w_i$ for $i=1,2,3,4$ be a basis of $H^1_{dR}(U,\nabla)$.
A period is the pairing of a cycle in $H^{irreg}_1(U,\nabla^*)$ and a cocycle
in $H^1_{dR}(U,\nabla)$,  which is given by the integration
$$
< \eta_i\otimes \exp y , \omega_j > := \int_{\eta_i} \exp y\cdot \omega_j. 
$$

The period determinant is the determinant of the period matrix:
$$
per(U,\nabla) := \det \begin{pmatrix}
\int_{\eta_i}\exp y \cdot \omega_j
\end{pmatrix}_{i,j=1,2,3,4}.
$$

Let $k$ be the field of definition of $U$. 
It is not a well-defined number in $\C$ since we can change the
basis of the $k$-vector space $H^1_{dR}(U,\nabla)$ by a matrix in $\GL(H^1_{dR}(U,\nabla))=\GL_4(k)$.
Moreover, we don't have a canonical choice of a basis in $H^{irreg}_1(U,\nabla^*)$.
It depends on the choice of local solution of $\nabla^*$. 

Throughout this article, we
choose $\exp y$ for the basis. Hence the period determinant is well-defined in
$\C^*/k^*$.

\section{Direct image connection}
We will keep the same convention as in the previous section.

In the previous section, we have expressed the period as an exponential elliptic
integration. As it is known, elliptic integrations are difficult to evaluate, thus
we take its direct image $\pi_*(\nabla)$ onto $\A^1$ where 
$\pi: (x,y) \mapsto y$ for $(x,y)\in U$.
Note $\nabla = \pi^*(d+dy)$. Using the projection formula, we have
\begin{equation}
\begin{split}
\pi_*(\sO_U, \nabla) =& \pi_* \pi^*(\sO_{\A^1}, d+dy) \\
	=& (\sO_{\A^1}, d+dy) \otimes \pi_*(\sO_U, d)
\end{split}
\end{equation}

Let $f(x)$ be $x(x-1)(x-\lambda)$ throughout this article
and $L$ be the function field of the elliptic curve
$k(x,y)/(y^2 - f(x))$. 
The trace map of $L$ in $k(y)$ induces a map
$$\pi_*(\sO_U) \stackrel{\Tr_{L/k(y)}}{\longrightarrow} \sO_{\A^1}.$$
Since this map has a section and flat, we have splitting of 
$(\pi_*(\sO_U), d)$ as connection: 
$(\pi_*(\sO_U),d)=(\sO_{\A^1}, d) \oplus (V, \nabla')$ for a rank $2$ connection
on the kernel of $\Tr_{L/k(y)}$, which is 
generated by $v_1 :=  - \frac{\lambda+1}3 + x$
and $v_2 := -\frac{\lambda^2 +1}3 + x^2$.
Therefore we have
$$
\pi_*(\sO_U,\nabla) = (\sO_{\A^1}, d+dy) \oplus ((\sO_{\A^1},d+dy) \otimes ( V, \nabla')).
$$
Let $x_1, x_2$ be the roots of $f'(x)$.
Note that $\pi_*(\sO_U,\nabla)$ obtains regular singularities via ramifications
at the points $(x_i, \pm\sqrt{f(x_i)})$ for $i=1,2$, which are concentrated
in $(V,\nabla')$.  

With respect to the basis $\{ v_1, v_2\}$, the connection matrix of $\nabla'$ is
$$
\frac{2ydy}{3(y^2 - f(x_1))(y^2 - f(x_2))}
\cdot 
\begin{pmatrix}
\frac19(\lambda+1)(2\lambda^2 - 3\lambda +2) + y^2 &
	\frac29 \lambda(\lambda+1) - \frac23(\lambda+1)y^2 \\
-\frac29(\lambda^2 -\lambda +1) & -\frac29 \lambda(\lambda+1) + 2y^2
\end{pmatrix}
$$
which has simple poles at the zeroes of $y= \pm\sqrt{f(x_i)}$ for $i=1,2$.

We have two different configurations of singularities after $\lambda$: one is
the generic case when $\lambda^2 - \lambda +1 \ne 0$ and the other is when
$\lambda^2 - \lambda +1= 0$. In the first case, we have four different singular points
at $y= \pm \sqrt{f(x_1)}, \pm \sqrt{f(x_2)}$. In the second case, we have only two singular points $y = \pm \sqrt{f(x_1)}$.

In each case, we will approximate the answer by a sequence of regular singular connections over $\P^1$.

\section{Product formula:Terasoma's work}

In this section, we recall the main theorem in \cite{te1}. The theorem 
tells the exact value of the period determinant of a regular singular connection on $\A^1$ for a  canonical choice of basis of de Rham cohomology and relative homology valued in the dual local system, assuming some extra conditions. This result will be applied 
for the approximation later.

We will firstly recall some necessary notions for the result.

Let $D = \{ \lambda_1, \lambda_2, \ldots, \lambda_n\} $ be $n$ distinct points in 
$\A^1$.  A logarithmic connection with poles at $D \cup \{\infty\}$ 
$$
\nabla : \sV \to \sV \otimes \Omega_{\P^1} (\log (D + \infty))
$$
on a trivial vector bundle $\sV$ of rank $r$ can be written as
$$
\nabla = d + \sum^n_{i=1} \frac{B^{(i)}}{x-\lambda_i},
$$
where $B^{(i)} = \Res_{\lambda_i}(\nabla) \in {\rm End}(\C^r)$ is the residue of $\nabla$ at $\lambda_i$. 
And  the residue at $\infty$ is 
$B^{(\infty)} := \sum^n_{i=1} (-B^{(i)})$. Throughout this section, 
we assume that no two eigenvalues of the residue are different by an integer and
that they have positive real part.

\begin{definition}
Let $(V,\nabla)$ be a connection on $\A^1 - D$.
$(\tilde{V},\tilde{\nabla})$ be its logarithmic extension to $\P^1$. If no eigenvalue of the residue of 
$\nabla$ at $D$ is a non-positive integer, it  is called a small extension of $\nabla$ along $D$.
\end{definition}

When $\nabla$ is small at $D$, then the de Rham cohomology $H^1_{dR}(\A^1, \nabla)$
of the logarithmic de Rham complex is generated by
$$
(\frac1{x-\lambda_i} - \frac1{x-\lambda_{i+1}}) dx \otimes w, 
\quad\text{for $i=1,2,\ldots,n-1$}
$$
where $w$ is a vector in $V$.
Thus $\rank H^1_{dR}(\A^1, \nabla)= r(n-1)$, which is canonically isomorphic to 
$H^1(\A^1, j_!j^*(\nabla))$, where $j:\A^1 - D \hookrightarrow \A^1$ is 
the open embedding (cf. \cite{st}, \cite{te1}). 

Now we need the corresponding homology theory to yield a perfect pairing with 
the previously described de Rham cohomology. This is the relative homology of the
pair $(\A^1, D)$ valued in the local system of $\nabla^*$. It is isomorphic to 
the relative homology of the pair $(U, \cup_i D_i)$ valued in the dual local system,
where $D_i$ is a sufficiently small disk around $\lambda_i$. 
It can be seen as Borel-Moore homology valued in the local system.
A cycle in this homology group can be written as linear sum of 
$\gamma_\alpha\otimes f_\alpha$, where $\gamma_\alpha$ is a cycle 
in $H_{i}(\A^1,D,\Z)$ and $f_\alpha$ is a branch of a solution of $\nabla^*$ 
over $\gamma_\alpha$. 
The main theorem tells us that $H_1(\A^1,D,\nabla^*)$ and $H^1_{dR}(\A^1,\nabla)$
make a perfect pair via integration. 

For a form  $\omega : = 
( \frac1{x-\lambda_i} - \frac1{x-\lambda_{i+1}})dx \otimes v$ in $H^1_{dR}(\A^1, \nabla)$,
and $\delta := \gamma\otimes f$ in $H_1(\A^1, D, \nabla^*)$, we define their pairing
as
$$
<\omega,\delta> := \int_\gamma(\frac1{x-\lambda_i} - \frac1{x-\lambda_{i+1}})<v,f> dx.
$$
 
The product formula \`a la Terasoma shows that the above pairing is perfect. 
To state the formula we need to introduce the tame symbol and the Gamma factor 
of a connection. 

For a rank $1$ connection $\nabla = d + \sum^n_{i=1} \frac{b^i}{x-\lambda_i}$ over
$\A^1$ with log poles at $\lambda_1,\ldots,\lambda_n$, where $b^i \in \C$ is the residue
at $\lambda_i$. Assume, as before, $b_i$ have positive real part. Let $p$ be a fixed point in $\A^1 - \{\lambda_1, \ldots, \lambda_n\}$ and $\gamma_i$ be a fixed path
from $p$ to $\lambda_i$. 
$\nabla^*$ has a multi-valued solution $\prod^n_{i=1} (x-\lambda_i)^{b^i}$. 
Fix a branch of $\log(x-\lambda_i)$ around $\lambda_i$ to 
have real value at $\lambda_i + \epsilon$ for small $\epsilon \in \R_{>0}$, thus we have
fixed the branch of $(x-\lambda_i)^{b^i}= \exp(b^i \log(x-\lambda_i))$ as well.
Let $D(x)_{\gamma_i}$ be the brach of the above solution of $\nabla^*$ 
on $\gamma_i$ according to the chosen branch of the logarithm.

The \textit{tame symbol} of the rank $1$ connection $\nabla$ is a value 
depending on the path from $p$ to $\lambda_i$
$$
(\nabla,(x-\lambda_i))_{\gamma_i} := \lim_{\begin{smallmatrix}x\to \lambda_i \\ 
	\text{along $\gamma_i$}\end{smallmatrix}} 
	\frac{D(x)_{\gamma_i}}{(x-\lambda_i)^{b^i}} = \prod_{j\ne i} (\lambda_i - \lambda)^{b^j}
$$
and for a path $\gamma_\infty$ from $p$ to $\infty$, 
$$
(\nabla, \frac1x)_{\gamma_\infty} := \lim_{x\to\infty} D(x)_{\gamma_\infty} x^{b^\infty}.
$$

For a higher rank connection, we define the \textit{tame symbol} as that of the
determinant connection.

Keeping the assumption that the residue $\lambda_i$ has only positive real part, 
we define the Gamma factor as a generalization of the Gamma function for the
residue:
$$
\Gamma_\lambda(\nabla) = 
\begin{cases}
\det (\int_0^\infty x^{\Res_\lambda(\nabla)} e^{-x} \frac{dx}x& \text{for $\lambda \in D$}, \\
\det (\int_0^\infty x^{-\Res_\infty(\nabla)} e^{-x} \frac{dx}x & \text{for $\lambda= \infty$}\\
\end{cases}
$$

Let $\nabla$ be a connection with log poles at $D=\{\lambda_1,\ldots,\lambda_n\}$ 
and $\infty$. 
Assume the eigenvalues of $\Res_{\lambda_i}$ has positive real parts and they are
not different by an integer. Let $\{e_i\}_{i=1,\ldots,r}$ be a basis for the underlying
vector space where the connection values and $\{ e^*_j\}$ be the dual basis.
Let $\delta_j := \gamma_j - \gamma_{j-1}$ and $\Sol(\nabla^*)(e^*_q)$ be the
branch of $\Sol(\nabla^*)$ analytically continued along $\gamma_i$ with the initial value
$e^*_q$ at $p$. We take basis
$$
\omega_i(e_p) := (\frac1{x-\lambda_i} - \frac1{x-\lambda_{i+1}})dx \otimes e_p
\quad \text{ in $H^1_{dR}(\A^1,\nabla)$,}
$$
and
$$
\delta_j(e^*_q) := \delta_j \otimes \Sol(\nabla^*)(e^*_q),
$$
to define a period matrix for $(i,j)$ as
$$
A_{ij} := \big( <\delta_j(e^*_q), \omega_i(e_p)>\big)_{1\le p,q \le r}.
$$

\begin{theorem}[Terasoma\cite{te1}]
\label{thm:period}
The determinant of the period matrix with respect to the above basis of 
$H^1_{dR}(\A^1, \nabla)$ and $H_1(\A^1,D,\nabla^*)$ is
$$
\det (A_{ij})_{1\le i,j\le n-1} = \prod^n_{i=1} (\nabla, x-\lambda_i)_{\gamma_i}\cdot
	(\nabla,\frac1x)^{-1}_{\gamma_\infty}\prod^n_{i=1} \Gamma_{\lambda_i}
	(\nabla)\cdot\Gamma_\infty(\nabla)^{-1}
$$
in $\C^*$. Therefore $H^1_{dR}(\A^1,\nabla)$ and $H_1(\A^1,D,\nabla^*)$ make a
perfect pairing because tame symbols and Gamma factors never vanish.
\end{theorem}
\begin{proof}
See \cite{te1}.
\end{proof}

We can take a different basis for $H^1_{dR}(\A^1, \nabla)$ to evaluate the period
determinant:
$$
\eta_i(e_p) := \frac{x^{i-1}}{\prod^n_{k=1} (x-\lambda_k)} \otimes e_p,
\quad \text{for $i=1,\ldots,n-1$ and $p= 1,\ldots, r$.}
$$
\begin{proposition}
\label{prop:per-det}
The period determinant of the pairing of $H^1_{dR}(\A^1,\nabla)$ and 
$H_1(\A^1, D, \nabla^*)$ with respect to the basis $\eta_i(e_p)$ and 
$\delta_j(e_q^*)$ is
$$
\det \big(  \int_{\delta_j(e^*_q)}( \eta_i(e_p) )\big) = 
\frac{\det(\int_{\delta_j(e_q^*)} \omega(e_p))}{\Delta(\lambda_1,\ldots,\lambda_n)^r}
$$
in $\C^*$ where 
$\Delta(\lambda_1,\ldots,\lambda_n) = \prod_{i<j}(\lambda_i - \lambda_j)$ 
the
vandermonde determinant of $(\lambda_1,\ldots,\lambda_n)$.
\end{proposition}
\begin{proof}
We will prove it for rank $1$ case. It will be directly generalize for higher rank cases.
Set $\eta'_i = \prod^{i-1}_{k=1} (x-\lambda_k) \eta_i$ for 
$\eta_i = \frac{x^{i-1}dx}{\prod^n_{k=1} (x-\lambda_i)}$.
Then we have
$
x\eta'_i = \eta'_{i+1} + \lambda_i \eta_i. 
$
By induction on $i$, it follows that $\eta_i = \eta_i' + \sum_{k>i} c_k \eta_k'$ for some constants
$c_k$. Hence
$$
\det(\int_{\delta_j} \eta_i ) = \det(\int_{\delta_j} \eta'_i).
$$
Whereas, for $a_i = \Res_{\lambda_i}  \eta'_i$,
$\eta'_i = \sum_{k=i}^n \frac{a_k dx}{x-\lambda_k} =  
\sum_{k=i}^{n-1} \sum_{\ell =i}^k a_\ell \omega_k + 
\sum_{\ell =i}^n \frac{a_\ell dx}{x-\lambda_n}.$ 
Since $\Res_{\infty}\eta_i' =0$,  $\sum_{k =i} a_i$ the
sum of all residues of $\eta'$ is $0$. 

$a_i = \Res_{\lambda_i} \eta'_i= \frac1{\prod_{k = i}^n (\lambda_i - \lambda_k)}$ and 
we have
$\eta'_i =  \frac{\omega_i}{\prod_{k \ge i} (\lambda_i - \lambda_k)} + 
\sum_{k > i} c_k \omega_k$ for some constants $c_k$.

Therefore the period determinant is
$$
\det(\int_{\delta_j} \eta_i ) = \frac1{\prod_{i<j}(\lambda_i -\lambda_j) }
\det(\int_{\delta_j} \omega_i),
$$
which finishes the proof for rank $1$ case. In higher rank case, we have multiple
contribution of the constant to its rank. Hence we obtain the power of the vandermonde matrix in the denominator.
\end{proof}

In \cite{te0} and \cite{va}, the authors treat the determinant of the period matrix 
$$
\bigg(
\int_{\delta_j} \prod^n_{k=1} (x-\lambda_k)^{s_k-1} x^{i-1} dx
\bigg)_{1\le i,j \le n-1},
$$
where $\delta_i$ is as above.
The integration of each entry is, in fact, the pairing
$$
\int_{\delta_j} \prod^n_{k=1} (x-\lambda_k)^{s_k-1} x^{i-1} dx =
< \delta_j \otimes \prod^n_{k=1}(x-\lambda_k)^{s_k} ,  \eta_i >
$$
for a rank $1$ connection $\nabla = d + \sum^n_{k=1} \frac{s_k}{x-\lambda_k}dx$ 
on $\A^1$ and $\eta_i$ are the same as we have described before.

Using Proposition \ref{prop:per-det}, we obtain the period determinant of the above
connection as
$$
\det\big( \int_{\delta_j} \prod^n_{k=1} (x-\lambda_k)^{s_k-1} x^{i-1} dx\big) =
\frac{\Gamma(s_1) \cdots \Gamma(s_n)}{\Gamma(s_1 + \ldots + s_n)}
\prod_{i<j}(\lambda_j - \lambda_i) \prod^n_{i=1}
(\prod_{j\ne i}(\lambda_j - \lambda_i)^{s_i-1}
$$
in $\C^*$.

\section{Period integral}
In this section, we want to evaluate the integration via an approximation.
As we have seen earlier, $\pi_*(\nabla)$ splits into 
$(\sO_{\A^1},d+dy)\oplus (V,\nabla')$ for a rank $2$ connection $\nabla'$.
Firstly, before applying the results of Terasoma for the approximation, 
we have to check if the rank $2$ part of $\pi_*(\nabla)$ satisfies the 
condition on eigenvalues of the residues.
Let $D=\{z_1, \cdots,z_n\}$ be a nonempty set of points in $\A^1$. For
a connection $\nabla$ over $\A^1$, possibly singular, we denote by $d_D$ the
standard exterior differentiation on $\sO_{\A^1}(-D)$ and define
$\nabla_D:= \nabla\otimes d_D$.

When $\lambda^2-\lambda+1 \ne 0$, the residue of $\nabla'$ at 
$\{\pm \sqrt{f(x_i)}\}$  has two eigenvalues $0$ and $-1/2$. This means
Terasoma's result is not directly applicable. In case $\lambda^2 -\lambda+1=0$,
the eigenvalues are $1/3, 2/3$. 

For the generic case $\lambda^2 -\lambda +1 \ne 0$, we modify $\nabla$ to
satisfy the eigenvalue condition, tensoring the rank $1$ connection 
$(\sO_{\A_1}(-D), d+\varpi)$ where 
$\varpi = (\frac1{y-\sqrt{f(x_1)}} + \frac1{y+\sqrt{f(x_1)}} + \frac1{y-\sqrt{f(x_2)}} +\frac1{y+\sqrt{f(x_2)}})dy$.
Then the new connection  obtained is
$(\sO_{\A^1}(-D), d + \varpi) \otimes \pi_*(\sO_U, \nabla) = 
\pi_*\pi^*(\sO_{\A^1}(-D), d+dy+\varpi)$ and this has regular
singularities at $D$ with the eigenvalues $1, 3/2$  of the residues at $D$.

We will calculate the period determinant for the above connection
$\pi_*(\nabla_D)$. The generic case $\lambda^2-\lambda+1\ne 0$
will be treated first. Since $\pi_*(\nabla_D)$ splits into a 
direct sum of a rank $1$ connection and a rank $2$ connection,
the period determinant is the product of the period of rank $1$ part
and that of rank $2$ part. The period determinant of a rank $1$ 
irregular connection of the form $\nabla =d+ dF +
\sum_{i=1}^n\frac{s_i dx}{x-\lambda_i}$ was evaluated by Terasoma:

\begin{theorem}[Terasoma\cite{te2}]
Let $\nabla = d + dF+ \sum_{i=1}^n\frac{s_i dx}{x-\lambda_i}$ be a connection
on $\A^1$ for a polynomial $F$. Let $I_i$ be a fixed path from $\lambda_i$ to $\lambda_{i+1}$ for $i=1,\ldots, n-1$ and $J_i$ be a fixed 
path from $\lambda_n$ to $\infty$ along the $i$-th rapid decay sector for $i=1,\ldots,
\deg (F)$.
Taking the basis $\eta_j = \frac{x^{j-1}dx}{\prod_{i=1}^n (x-\lambda_i)}$
of $H^1_{dR}(\A^1, \nabla)$ and the basis 
$\delta_j :=I_j\otimes \prod_{k=1}^n (x-\lambda_k)^{s_k} \exp (F(x))$ ($i=1,\ldots,n-1$) and
$J_j \otimes \prod_{k=1}^n (x-\lambda_k)^{s_k} \exp (F(x))$ 
($j=1,\ldots,\deg (F)$), we have the determinant of the period matrix
\begin{equation}
\begin{split}
D =& \det ( <\delta_j, \eta_i> )\\
	=& (2\pi)^{(d-1)/2}\Gamma(s_1)\cdots\Gamma(s_n)(da_d)^{-s-(d-1)/2}
		(-1)^{ds+d(d-1)/4}\\
		&\times \prod^n_{i=1}(\prod_{j\ne i} 
		(\lambda_i-\lambda_j)^{s_i-1}\prod_{i<j}(\lambda_j-\lambda_i)) \\
		&\times \prod^n_{i=1} \exp(F(\lambda_i)) \prod_{F'(u)=0} \exp(F(u)),
\end{split}
\end{equation}
in $\C^*$, where $s=s_1 +\cdots + s_n$.
\end{theorem}

The rank $1$ part of $\pi_*(\nabla_D)$ is 
$d+ dy + \varpi$. 
Applying the formula, we have
the period determinant of the rank $1$ part
\begin{equation}
\begin{split}
D =& (\sqrt{f(x_1)} + \sqrt{f(x_2)})(2\sqrt{f(x_1)}2\sqrt{f(x_2)})(\sqrt{f(x_1)}-\sqrt{(x_2)})\\
& \times (\sqrt{f(x_1)} - \sqrt{f(x_2)})(-2\sqrt{f(x_2)})(-\sqrt{f(x_2)}-\sqrt{f(x_1)}) \\
= & 2^2 \sqrt{f(x_1)}\sqrt{f(x_2)} (f(x_1)-f(x_2))^2
\end{split}
\end{equation}
in $\C^*$.

For simplicity, we will denote 
$$
\omega = \frac{2ydy}{3(y^2 -f(x_1))(y^2 -f(x_2))}\in \Gamma(\A^1,\Omega(\log D)).
$$
The rank $2$ part has the connection
\begin{equation}
\begin{split}
\nabla_1 = & (d+ dy)\otimes \nabla_D' \\ =& d + I dy + \varpi \\
	& + 
	\begin{pmatrix}
		\frac19(\lambda+1)(2\lambda^2 - 3\lambda +2)+ y^2 &
		\frac29\lambda(\lambda^2 +1) -\frac23 (\lambda+1) y^2 \\
		-\frac29(\lambda^2 - \lambda+1) & 
		-\frac29\lambda(\lambda+1) + 2y^2
	\end{pmatrix} \omega.
\end{split}
\end{equation}

We are now prepared to approximate the integrations. 
Note that the residue of the new rank $2$ connection has $1,3/2$ at 
each point of $D$. 

Let $e_1,e_2$ be the standard basis of $\Gamma(\P^1,\sO^{\oplus 2})$ and $e^*_1,e^*_2$
be the dual so that $<e_i,e_j^*> = \delta_{ij}$. Let $\gamma_i$ be a fixed
path from $0$ to $-\sqrt{f(x_1)}, \sqrt{f(x_1)}, -\sqrt{f(x_2)}, \sqrt{f(x_2)}$ respectively,
for $i=1,2,3,4$. $\gamma_\infty$ is a path from $0$ to $\infty$ along the rapid
decaying sector around $\infty$.  $I_j$ is the chain $\gamma_{i+1} - \gamma_i$
for $i=1,2,3$ and $I_4 := \gamma_\infty - \gamma_4$.
With this notation, we have the following cycles as a basis for 
$H^{irreg}_1(\A^1,D, \nabla_1^*)$:
$$
C_{j,b} := I_j \otimes (\exp y\cdot \Sol(\nabla^*_1))(e^*_b), \quad
\text{for $j=1,\ldots,4$ and $b=1,2$}
$$
and $H^1_{dR}(\A^1, \nabla_1)$ has a basis
$$
\eta_{i,a} := \frac{y^{i-1}dy}{(y-\sqrt{f(x_1)})(y+\sqrt{f(x_1)})(y-\sqrt{f(x_2)})
	(y+\sqrt{f(x_2)})}
	\otimes e_a
$$
for $i=1,2,3,4$ and $a=1,2$.

The pairing of $\eta_{i,a}$ and $C_{j,b}$ is the integration
\begin{equation}
\begin{split}
P_{i,j,a,b} :=& <\eta_{i,a},C_{j,b}> \\
= &\int_{I_j} \frac{y^{i-1} \exp y <e_a, \Sol(\nabla^*_1)(e_b^*) > dy}
	{ (y-\sqrt{f(x_1)})(y+\sqrt{f(x_1)})(y-\sqrt{f(x_2)})(y+\sqrt{f(x_2)})}.
\end{split}
\end{equation}

Let $I_j^{(m)}$ be $I_j$ for $i=1,2,3$ and $I_4^{(m)}$ be $\gamma_{(-m)}-\gamma_4$,
where $\gamma_{(-m)}$ is a path from $0$ to $-m$ sitting in the rapid decay sector
of $\nabla_1^*$ near $\infty$ and converges uniformly to $\gamma_\infty$. 
The above integration is approximated by 
\begin{equation}
\begin{split}
\mbox{}&P_{(m),i,j,a,b}\\
=& \int_{I_j^{(m)}} \frac{y^{i-1}(1+\frac{y}m)^m < e_a, \Sol(\nabla^*_1) (e^*_b) > dy}
	{(y-\sqrt{f(x_1)})(y+\sqrt{f(x_1)})(y-\sqrt{f(x_2)})(y+\sqrt{f(x_2)})} \\
= & (\frac1m)^m \int_{I_j^{(m)}} 
	\frac{y^{i-1}(y+m)^{m+1} <e_a,\Sol(\nabla_1^*)(e^*_b) > dy}
		{(y-\sqrt{f(x_1)})(y+\sqrt{f(x_1)})(y-\sqrt{f(x_2)})(y+\sqrt{f(x_2)})}
\end{split}
\end{equation}
as $m$ tends to $\infty$.

This appears as the period integration of a regular singular connection
\begin{equation}
\begin{split}
\nabla_{(m)} =& (d+  \frac{(m+1)dy}{y+m}) \otimes \nabla'   \\
	 =& d + I\frac{(m+1)dy}{y+m} + I\varpi \\
	& + 
\begin{pmatrix}
		\frac19(\lambda+1)(2\lambda^2 - 3\lambda +2)+ y^2 &
		\frac29\lambda(\lambda^2 +1) -\frac23 (\lambda+1) y^2 \\
		-\frac29(\lambda^2 - \lambda+1) & 
		-\frac29\lambda(\lambda+1) + 2y^2
	\end{pmatrix} \omega.
\end{split}
\end{equation}

Note that $\frac{(m+1)dy}{y+m}$ was chosen  to approximate $dy$.
For each $m$, calculating the period determinant of $\nabla_{(m)}$,
we approximate that of $\nabla_1$. We obtain the Tame symbols 
of $\nabla_{(m)}$ as follows:
\begin{equation}
\begin{split}
(\nabla_{(m)}, y-\sqrt{f(x_1)})_{\gamma_1} &= (2\sqrt{f(x_1)})^{5/2}
	(f(x_1)-f(x_2))^{5/2} (\sqrt{f(x_1)} + m)^{2(m+1)},\\
(\nabla_{(m)}, y+\sqrt{f(x_1)})_{\gamma_2} &=
	(-2\sqrt{f(x_1)})^{5/2} (f(x_1)-f(x_2))^{5/2}(-\sqrt{f(x_1)} + m)^{2(m+1)},\\
(\nabla_{(m)}, y-\sqrt{f(x_2)})_{\gamma_3} &= (2\sqrt{f(x_2)})^{5/2}
	(f(x_2)-f(x_1))^{5/2} (\sqrt{f(x_2)} + m)^{2(m+1)},\\
(\nabla_{(m)}, y+\sqrt{f(x_2)})_{\gamma_4} &=
	(-2\sqrt{f(x_2)})^{5/2} (f(x_2)-f(x_1))^{5/2}(-\sqrt{f(x_2)} + m)^{2(m+1)},\\
(\nabla_{(m)}, y+m)_{\gamma_{-m}} &= (m^2 - f(x_1))^{5/2}
	(m^2 -f(x_2))^{5/2} \quad \text{and}\\
(\nabla_{(m)}, \frac1y)_{\gamma_\infty} &= 1.
\end{split}
\end{equation}

The Gamma factors are
$$
\Gamma_{\pm \sqrt{f(x_i)}}(\nabla_{(m)}) = \Gamma(1)\Gamma(\frac32) = \frac12\Gamma(\frac12)
=\frac{\sqrt{\pi}}2,
$$
$$
\Gamma_{-m}(\nabla_{(m)}) = \Gamma(m+1)^2
$$
and
$$\Gamma_\infty (\nabla_{(m)}) = \Gamma((m+1)+\frac{14}3)
	\cdot\Gamma((m+1)+\frac{16}3).
$$

With the above factors, we have the period determinant
for the regular singular connection $\nabla_{(m)}$ with
respect to the basis chosen below.
\begin{equation}
\begin{split}
C^{(m)}_{j,b} : =& I^{(m)}_j\otimes 
\Sol(\nabla_{(m)}^*) (e^*_b)\\
\omega^{(m)}_{1a} :=& (\frac1{y+\sqrt{f(x_1)}} - \frac1{y-\sqrt{f(x_1)}})dy\otimes e_a,\\
\omega^{(m)}_{2a} :=&  (\frac1{y-\sqrt{f(x_1)}} - \frac1{y+\sqrt{f(x_2)}})dy\otimes e_a,\\
\omega^{(m)}_{3a} := & (\frac1{y+\sqrt{f(x_2)}} - \frac1{y-\sqrt{f(x_2)}})dy\otimes e_a,\\
\omega^{(m)}_{4a} := & (\frac1{y-\sqrt{f(x_2)}}-\frac1{y+ m} ) dy\otimes e_a.
\end{split}
\end{equation}

The period determinant of $\nabla_{(m)}$ with respect to $\omega^{(m)}_{i,a}$ 
and $C^{(m)}_{j,b}$ obtained by the application of Theorem \ref{thm:period}
as follows:
\begin{equation}
\begin{split}
D_{(m)} =& \det ( < C^{(m)}_{j,b}, \omega^{(m)}_{i,a}> ) \\ 
	=& 2^6 \pi^2 f(x_1)^{5/2} f(x_2)^{5/2} (f(x_1)-f(x_2))\\
	& \times (m^2 - f(x_1))^{2(m+1)+5/2} (m^2 - f(x_2))^{2(m+1)+5/2}\\
	& \times \frac{\Gamma(m+1)^2}{\Gamma((m+1)+\frac{14}{3})
		\Gamma((m+1)+\frac{16}3)}.
\end{split}
\end{equation}

The above value will not converge, but
using Proposition \ref{prop:per-det} we obtain a convergent sequence of
period determinants
\begin{equation}
\begin{split}
P_{(m)} := 
\det (P_{(m),i,j,a,b}) = (\frac1m)^{8m} \frac{D_{(m)}}{\Delta^2_{(m)}}
\end{split}
\end{equation}
where $\Delta_{(m)}$ is the Vandermonde determinant for $\pm\sqrt{f(x_1)}$,
$\pm\sqrt{f(x_2)}$, $-m$ and 
$\Delta_{(m)}$ is equal to
$$
- 2^2 \sqrt{f(x_1)}\sqrt{f(x_2)} (f(x_2)-f(x_1))^2 (m^2 - f(x_1))(m^2 - f(x_2)).
$$

Together with the above, the period determinant of $\nabla_{(m)}$ with respect to
$\eta^{(m)}_{i,a}$ and $C^{(m)}_{j,b}$ is
\begin{equation}
\begin{split}
P_{(m)} =& \det (P_{(m),i,j,a,b}) = \det (<\eta^{(m)}_{i,a}, C^{(m)}_{j,b}>) \\
	=& (\frac1m)^{8m} \frac{D_{(m)}}{\Delta^2_{(m)}} \\
	=& \frac{2^2 \pi^2 f(x_1)^{3/2} f(x_2)^{3/2}}{ (f(x_1) -f(x_2))^{3}}\\
	&\times (\frac1m)^{8m} (m^2 - f(x_1))^{2m}(m^2 -f(x_2))^{2m} \\
	&\times m^{-10} (m^2 - f(x_1))^{5/2} (m^2 -f(x_2))^{5/2} \\
	& \times m^{10} \frac{\Gamma(m+1)^2}{\Gamma((m+1)+\frac{14}3)
		\Gamma((m+1)+\frac{16}3)},
\end{split}
\end{equation}
which converges to $P$ the period determinant of $\nabla_1$.
Each term in the above converges as follows:
\begin{equation}
\begin{split}
(\frac{1}{m})^{8m} (m^2 - f(x_1))^{2m}(m^2-f(x_2))^{2m} &\to 1\\
m^{-10} (m^2 - f(x_1))^{5/2}(m^2-f(x_2))^{5/2} & \to 1\\
m^{10}\frac{\Gamma(m+1)^2}{\Gamma((m+1)+\frac{14}3)\Gamma((m+1)+\frac{16}3)} &\to 1.
\end{split}
\end{equation}

Therefore  $P$ is
$$
\frac{2^2\pi^2 f(x_1)^{3/2}f(x_2)^{3/2}}{(f(x_1)-f(x_2))^3}.
$$
Note the above value is found in $\C^*$.

The period determinant of $\pi_*(\nabla)_D = \pi_*\pi^*(d+dy+\varpi)$ is then,
\begin{equation}
\begin{split}
&(\text{Period of $(\sO_{\A^1}, d+ dy +\varpi)$})\times
(\text{Period of $\nabla_1$})\\
=& \frac{2^4 \pi^2 f(x_1)^2 f(x_2)^2}{f(x_1)-f(x_2)}.
\end{split}
\end{equation}

\section{Comparision}
Recall that $U$ is the affine Legendre elliptic curve defined by the equation
$y^2 = x(x-1)(x-\lambda)$ and $\nabla = d+dy$. Let $D$ be as in the previous 
section and $\Sigma := \pi^{-1} D$ in $U$.
Then using the projection formula, 
$\pi_*(\nabla)_D = d_D \otimes \pi_*\nabla$, where $d_D$ and $\nabla_D$ are 
defined as before. In the same manner, we denote
by $\nabla_\Sigma$, the twisted connection  $d_\Sigma\otimes \nabla$ on $\sO(-\Sigma)$.

Then $\pi$ induces a canonical isomorphism of de Rham cohomologies:
$$
\pi^*: H^1_{dR}(\A^1, \nabla_D) \to H^1_{dR}(U,\nabla_\Sigma).
$$
If $\omega = \pi^*\eta$ a de Rham form in $H^1_{dR}(U,\nabla_\Sigma)$
for a $\eta$ in $H^1_{dR}(\A^1,\nabla_D)$, then the functoriality of the pairing
implies $<\gamma, \omega> = <\gamma,\pi^*\eta> = <\pi_*\gamma,\eta>$.
Hence the period of $\pi_*(\nabla)_D$ that we calculated previously is
the same as that of $\nabla_\Sigma$ on $U$.

The following short exact sequence of de Rham complexes                         
\begin{equation}                                                                
\begin{CD}                                                                      
\sO_U(-\Sigma)  @>>>            \sO_U   @>>>     \sO_\Sigma \\                  
@VV{\nabla_\Sigma}V     @VV{\nabla}V                            @VVV \\         
\Omega_U(-\Sigma)(\log \Sigma) @=  \Omega_U     @>>>  0                         
\end{CD}                                                                        
\end{equation}                                                                  
yields 
\begin{equation}
H^0_{dR}(U,\nabla)\to H^0(U,\sO_\Sigma) \stackrel{\nabla}{\to} 
H^1_{dR}(U,\nabla_\Sigma) \to H^1_{dR}(U,\nabla) \to H^1 (U, \sO_\Sigma)
\end{equation}

Since $\nabla$ admits no (single-valued) solution on $U$, 
$H^0_{dR}(U,\nabla)=0$ and 
$H^1(U,\sO_\Sigma)$ vanishes by dimension reason, $H^1_{dR}(U,\nabla_\Sigma)$
is an extension of $H^1_{dR}(U,\nabla)$ by $H^0(U,\sO_\Sigma)$.

In the homology side, $H_1^{irreg}(U,\Sigma,\nabla^*)$ is isomorphic to
$H^{irreg}_1(\A^1,D,\pi_*(\nabla^*))$ via $\pi_*$.  The relative cohomology 
appears in the short exact sequence of the singular complexes valued in
$\nabla^*$:
\begin{equation}
0 \to \sC^{irreg}_* (\Sigma,\nabla^*) \to \sC^{irreg}_*(U,\nabla^*) \to
\sC^{irreg}_*(U,\Sigma,\nabla^*) \to 0,
\end{equation}  
which also yields a long exact sequence of homologies:
\begin{equation}
H^{irreg}_1(\Sigma,\nabla^*) \to H^{irreg}_1 (U,\nabla^*) \to 
H^{irreg}_1(U,\Sigma,\nabla^*) \to H^{irreg}_0 (\Sigma,\nabla^*) 
\to H^{irreg}_0(U,\nabla^*) 
\end{equation}

In the above sequence, as in cohomology side, we have vanishing
$H^{irreg}_1(\Sigma,\nabla^*)$ and $H^{irreg}_0(U,\nabla^*)$, respectively, by dimension reason and 
by the duality of de Rham cohomology and homology for irregular connections.

\begin{theorem}
Let $\gamma\otimes \Sol(\nabla^*)$ be a cycle in $H^{irreg}_1(U,\nabla^*)$,
thus a cycle in $H^{irreg}_1(U,\Sigma,\nabla^*)$. Suppose $\omega = \nabla f$
for a function $f$ in $H^0(U,\sO_\Sigma)$, then the pairing of 
$\gamma\times \Sol(\nabla^*)$ with $\nabla f$ is $0$.
\end{theorem}

\begin{proof}
The pairing is the integration
$<\gamma\otimes \Sol(\nabla^*),\nabla> = \int_\gamma <\Sol(\nabla^*),\nabla f>$.
This is equal to 
\begin{equation}
\begin{split}
&\int_{\gamma} (d <\Sol(\nabla^*),f> - <\nabla^*(\Sol(\nabla^*)),f>)\\
=&\int_\gamma d<\Sol(\nabla^*),f> = \int_{\partial\gamma} <\Sol(\nabla^*),f>
\end{split}
\end{equation}
by Stokes' theorem.

For a closed cycle $\gamma$, $\partial\gamma =0$. In this case we have nothing 
to prove. Otherwise, a tubular neighborhood of $\partial\gamma$ in $\gamma$ lies 
in the rapid decay sector of $\nabla^*$ around $\infty$, thus the integration 
yields $0$.
\end{proof}

Using the previous theorem, we conclude the following:
\begin{corollary}
The period determinant of $\nabla_\Sigma$ is
$$
per(U,\Sigma,\nabla_\Sigma) = per(U,\nabla) \times per(\Sigma,\nabla)
$$
in $\C^*/k^*$.
\end{corollary}

$\Sigma$ is given by the principal ideal $I:= ((y^2 -f(x_1))(y^2 -f(x_2)))$ 
in $R:= k[x,y]/(y^2 - f(x))$. It follows $H^0(U,\sO_\Sigma) = R/I$. So we
have a basis: $\{1,y$, $y^2,y^3$, $x,xy$, $xy^2, xy^3\}$ of $H^0(U,\sO_\Sigma)$,
whereas 
$H^{irreg}_0(\Sigma,\nabla^*)$ has a basis $\{p\otimes \exp y| p\in \Sigma\}$.

Let $x_3$ (resp. $x_4$) be the root of $f(x)-f(x_1)$ (resp. of $f(x)-f(x_2)$) such
that $x_3 \ne x_1$ (resp. $x_4 \ne x_2$). For a point $p = (x_i, \pm\sqrt{f(x_i)})$,
$<p\otimes\exp y, x^a y^b> = x^a_i (\pm\sqrt{f(x_i)})^b \cdot \exp (\pm\sqrt{f(x_i)})$,
for $a=0,1$ and $b=1,2,3,4$. $M_i$ is the following $(2\times 4)$-matrix:
$$
M_i := \begin{pmatrix}
1 & \sqrt{f(x_i)} & f(x_i) & f(x_i)\sqrt{f(x_i)} \\
1 & -\sqrt{f(x_i)} & f(x_i) & -f(x_i)\sqrt{f(x_i)}.
\end{pmatrix}.
$$

The period matrix in consideration is
$$
Q:= L \begin{pmatrix}
M_1 & x_1 M_1 \\
M_2 & x_2 M_2 \\
M_1 & x_3 M_1 \\
M_2 & x_4 M_2 
\end{pmatrix},
$$
where $L = (l_{ij})$ is a $(8,8)$-diagonal matrix with entries
$$
l_{ii} = \begin{cases}
\exp (f(x_1)) & \text{for $i\equiv 1 \pmod{4}$} \\
\exp (-f(x_1) ) & \text{for $i\equiv 2 \pmod{4}$} \\
\exp (f(x_2)) & \text{for $i\equiv 3 \pmod{4}$} \\
\exp (-f(x_2) ) & \text{for $i\equiv 0 \pmod{4}$} \\
\end{cases}
$$

It follows 
\begin{equation}
\begin{split}
\det Q &= \det
\begin{pmatrix}
M_1 & x_1 M_1 \\
M_2 & x_2 M_2 \\
M_1 & x_3 M_1 \\
M_2 & x_4 M_2 
\end{pmatrix} \\
&= \det\begin{pmatrix}
M_1 & x_1 M_1 \\
M_2 & x_2 M_2 \\
0 & (x_3 - x_1) M_1 \\
0 & (x_4 - x_2) M_2
\end{pmatrix} \\
& = (x_3 - x_1)^2 (x_4 -x_2)^2 \Delta^2_\Sigma,
\end{split}
\end{equation}
where $\Delta_\Sigma$ is the Vandermonde determinant for $\sqrt{f(x_1)}$,
$-\sqrt{f(x_1)}$, $\sqrt{f(x_2)}$ and $-\sqrt{f(x_2)}$.
This is $\Delta_\Sigma = 2^2 \sqrt{f(x_1)} \sqrt{f(x_2)} (f(x_2) - f(x_1))^2$.

$x_1, x_3$ are the two roots of 
$$
\frac{f(x)-f(x_1)}{x-x_1} = x^2 + (x_1 - (\lambda+1))x + x_1^2 - (\lambda+1) x_1 
+ \lambda.
$$
So,
$$
(x_1 - x_3)^2 = \frac{-3x_1^2 + 2(\lambda+1) x_1 + (\lambda+1)^2 }4.
$$

By the same way, we obtain
$$
(x_2-x_4)^2 = \frac{-3x_2^2 + 2(\lambda+1) x_2 + (\lambda+1)^2 }4
$$
and thus we have
\begin{equation}
(x_1 - x_3)^2 (x_2 - x_4)^2 = 2^{-4} (\lambda^2 - \lambda+ 1)^2.
\end{equation}

Altogether,
\begin{equation}
\begin{split}
\det Q &= (x_3 - x_1)^2 (x_4 - x_2)^2 \Delta_\Sigma^2 \\
	&= (\lambda^2 -\lambda + 1)^2 f(x_1) f(x_2) (f(x_2) - f(x_1))^4.
\end{split}
\end{equation}

After all, we obtain the period of $\nabla = d+dy$ over $U$:
$$
\frac{2^2 \pi^2 f(x_1) f(x_2)}{(\lambda^2 -\lambda +1)^2 (f(x_1) - f(x_2))^2}.
$$

To see the value explicitly, we need the values of $f(x_1) -f(x_2)$ and
$f(x_1)f(x_2)$.
$$
f(x_1) f(x_2) = -\frac1{3^3} \lambda^2 (\lambda-1)^2,
$$
and
\begin{equation*}
\begin{split}
(f(x_1) - f(x_2))^2 =& (f(x_1) + f(x_2))^2 - 4f(x_1)f(x_2)\\
	=& \frac{2^4}{3^6} (\lambda^2 -\lambda +1)^3.
\end{split}
\end{equation*}
It follows 
$f(x_1) -f(x_2) = \frac{2^2}{3^3} (\lambda^2 -\lambda +1)\sqrt{\lambda^2 -\lambda+1}$, for a suitable choice of the branch of the square root.

Finally, we have the period explicitly:
\begin{theorem}
The period determinant of $\nabla = d+dy$ over the affine Legendre elliptic curve
$y^2 = x(x-1)(x-\lambda)$ for $\lambda \ne 0,1, \frac{1\pm \sqrt{-3}}2$ is
$$
-2^{-6}3^{12}\pi^2 \frac{\lambda^2 (\lambda-1)^2}
{(\lambda^2 -\lambda+1)^9 \sqrt{\lambda^2 -\lambda +1}}
$$
in $\C^*/k^*$.
\end{theorem}

\section{Exceptional case: $\lambda^2 -\lambda+1 = 0$}

Finally, we handle the case $\lambda^2 - \lambda + 1 = 0$.
Recall that $x_1 = \frac{\lambda+1}3$ is a double root of $f'(x)$ as well
as a triple root of $f(x)-f(x_1)$.

Again using the projection formula, we see 
$\pi_*(\nabla) \simeq (\sO_{\A^1}, d+ dy) \oplus
(\ker \Tr, (d + dy)\otimes \nabla')$, 
where $\nabla'$ is, as before, the connection on the rank $2$ part of 
$\pi_*(\sO_U, d)$:
$$
\nabla'= d +
\begin{pmatrix}
1 & -\frac23 (\lambda +1) \\
0 & 2 
\end{pmatrix} \omega'
$$
with $\omega' = \frac{2ydy}{3(y^2 - f(x_1))}$.

$(d+dy)\otimes \nabla'$ has regular singularities at 
$D = \{\pm\sqrt{f(x_1)}\}$.
At each point of $D$, the residue of $\nabla'$ has two eigenvalues $1/3,2/3$.

Note $H^1_{dR}(\A^1, d+dy)$ is trivial, so the period of $\pi_*(\nabla)$ is
equal to that of $(\ker \Tr, (d+dy)\otimes \nabla')$. 

Let $\Sigma := \pi^{-1} D = \{ x_1, \pm\sqrt{f(x_1)}\} \subset U$.  
As we have seen in the previous section,
we have the following short exact sequence:
$$
0 \to H^0(\Sigma,\sO_\Sigma)\stackrel{\nabla}{\to}
H^1_{dR} (U,\nabla_\Sigma) \to H^1_{dR} (U,\nabla) \to 0.
$$
Note that $\pi^*$ induces a canonical isomorphism between
$H^1_{dR}(U,\nabla_\Sigma)$ and $H^1_{dR}(\A^1, \pi_*(\nabla)_D)$.

In the other hand the homologies make the short exact sequence:
$$
0 \to H^{irreg}_1 (U,\nabla^*) \to H^{irreg}_1(U,\Sigma, \nabla^*)
\stackrel{\delta}{\to} H^{irreg}_0 (\Sigma, \nabla^*) \to 0.
$$

Therefore the period in consideration $per(U,\nabla)$ satisfies
$$
per(U,\nabla_\Sigma) = per(\Sigma,\nabla)\cdot per(U,\nabla).
$$
Using the identification $ per(U,\nabla_\Sigma)= per(\A^1, \pi_*(\nabla)) $ will be approximated as in the previous section. 

For the approximation, we have to find  the cycles generating 
$H^{irreg}_1(\A^1,((d+dy)\otimes \nabla')^*)$.
Let $\gamma_1,\gamma_2$ be a fixed path from $0$ to $-\sqrt{f(x_1)}$
and to $\sqrt{f(x_1)}$  in $\A^1$ respectively. $\gamma_m$ is a path from
$0$ to $-m$. Define $I_1^{(m)} := \gamma_2 - \gamma_1$ and
$I_2^{(m)} := \gamma_m -\gamma_2$. Note that $I^{(m)}_1$ is independent of
$m$ and $I_2^{(m)}$ goes to $\infty$ along the rapid decay sector of $\nabla^*$.
Let $I_1 = I_1^{(m)}$ and $I_2 = \lim_{m\to \infty} I_2^{(m)}$.
Then $I_i\otimes \Sol(\nabla^*)(e_p^*)$ for $i,p=1,2$ make a basis for 
$H^{irreg}_1(\A^1, ((d+dy)\otimes \nabla')^*)$. 

In the other hand, $\eta_i(e_q):= \frac{y^{i-1}}{y^2 -f(x_1)} dy \otimes e_q$
for $i=1,2$, which is a basis for $H^1_{dR}(\A^1, (d+dy)\otimes\nabla')$. 

Let $\eta_i$
be $\frac{y^{i-1}dy}{y^2 -f(x_1)}$. We introduce a new connection 
$\nabla^{(m)} :=  (d+\frac{m+1}{m}dy)\otimes \nabla'$ which approximates
the rank $2$ part of $\pi_*(\nabla)$.
Then the $(i,j,p,q)$ entry of the period matrix is
\begin{equation}
\begin{split}
&<I_j \otimes \exp y \Sol(\nabla'^*)(e_p^*), \eta_i(e_q) > = 
\int_{I_j} <\Sol(\nabla'*_1) (e_p^*), e_q >  \exp y \cdot \eta_i \\
& = \lim_{m\to\infty} \frac1{m^m}\int_{I_j^{(m)}} <\Sol(\nabla'^*)(e_p^*),e_q>
	\frac{(y+m)^m y^{i-1}}{y^2 - f(x_1)} dy \\
& = \lim_{m\to \infty} \frac1{m^m} 
	<I_j^{(m)} \otimes \Sol (\nabla^{(m)})(e_q^*), \eta_i^{(m)} (e_p)> 
\end{split}
\end{equation}
where $\eta_i^{(m)}(e_p)$ is the de Rham form 
$\frac{y^{i-1}dy}{(y^2 - f(x_1))(y+m)}\otimes e_p$.

We denote the determinant of the above period matrix by $P_(m)$.

Applying the product formula, we obtain
\begin{equation}
\begin{split}
P_{(m)} =& \frac1{m^{4m}\Delta^2_{(m)}}
	(\nabla^{(m)}, y-\sqrt{f(x_1)})_{\gamma_1} \cdot
	(\nabla^{(m)}, y+\sqrt{f(x_1)})_{\gamma_2} \\
	&\times (\nabla^{(m)}, y+m)_{\gamma_m} \cdot (\nabla^{(m)}, \frac1y )^{-1}_\infty\\
	& \times \Gamma(\nabla^{(m)})_{-\sqrt{f(x_1)}} \cdot
		\Gamma(\nabla^{(m)})_{\sqrt{f(x_1)}} \cdot \Gamma(\nabla^{(m)})_{-m}
		\cdot \Gamma(\nabla^{(m)})^{-1}_\infty.
\end{split}
\end{equation}

The tame symbols are
\begin{equation}
(\nabla^{(m)}, y - \sqrt{f(x_1)})_{\gamma_1} =
	2\sqrt{f(x_1)} (\sqrt{f(x_1)}+ m)^{2(m+1)},
\end{equation}
\begin{equation}
(\nabla^{(m)},y + \sqrt{f(x_1)}_{\gamma_2}) =
	-2 \sqrt{f(x_1)} (- \sqrt{f(x_1)} + m)^{2(m+1)},
\end{equation}
\begin{equation}
(\nabla^{(m)}, y+m)_{\gamma_m} = m^2 - f(x_1),
\end{equation}
and
\begin{equation}
(\nabla^{(m)}, \frac1y)_\infty  = 1.
\end{equation}

And the Gamma factors are
\begin{equation}
\Gamma_{\pm\sqrt{f(x_1)}} (\nabla^{(m)}) = \Gamma(\frac13)\Gamma(\frac23),
\end{equation}
\begin{equation}
\Gamma_{-m} ( \nabla^{(m)}) = \Gamma(m+1)^2,
\end{equation}
\begin{equation}
\Gamma_\infty (\nabla^{(m)}) = \Gamma(\frac23 + (m+1)) \cdot 
			\Gamma(\frac43 +(m+1)).
\end{equation}

The period determinant $P_{(m)}$ is the product of the followings, which
converges:
\begin{equation}
\Gamma(\frac13)^2 \Gamma(\frac23)^2 
\end{equation}
\begin{equation}
\frac{(m^2 - f(x_1))^{2m} }{m^{4m}} \to 1
\end{equation}
\begin{equation}
\frac{m^2 - f(x_1)}{m^2} \to 1
\end{equation}
\begin{equation}
\frac{m^2 \Gamma(m+1)^2}{\Gamma(\frac23 + (m+1))\Gamma(\frac43+(m+1))}
\to 1
\end{equation}

From the functional equation $\Gamma(z)\Gamma(1-z) =\frac\pi{\sin \pi z}$ 
satisfied by the Gamma function, we obtain
$\Gamma(\frac13)^2 \Gamma(\frac23)^2  = (\pi/\sin(\frac\pi 3))^2 = 
\frac{2^2 \pi^2}3$.

Therefore, $per(U,\nabla_\Sigma) = \Gamma(\frac13)^2 \Gamma(\frac23)^2 = 
\frac{2^2\pi^2}3$.

Since $1,y$ is a basis of $\Gamma(U,\sO_\Sigma)$, 
$$
per(\Sigma,\nabla) = \det \begin{pmatrix} 1 & -\sqrt{f(x_1)} \\ 1 &\sqrt{f(x_1)}
\end{pmatrix} = 2\sqrt{f(x_1)}.
$$

From $x_1 = \frac{\lambda+1}3$ and $\lambda|^2 - \lambda+1 =0$, it follows
$f(x_1) = \frac{2\lambda-1}9 = \pm \frac{\sqrt{-3}}9$.

Finally, we obtain the period determinant
\begin{equation}
\begin{split}
per(U,\nabla) & = \frac{per(U,\nabla_\Sigma) }{per (\Sigma,\nabla)}
	= \frac{2^2 \cdot 3^{-1} \pi^2}{2\cdot 3 (-3)^{1/4}} = 
	\frac{2\pi^2}{(-3)^{1/4}}
\end{split}
\end{equation}
in $\C^*/ k^*$.

\end{document}